\input amstex
\documentstyle{amsppt}
\magnification=1200
\hcorrection{.25in}
\overfullrule=0pt

\redefine\R{{{\bold R}}}
\redefine\Z{{{\bold Z}}}

\define\hy{{\hbox{\it-}}}
\define\spec{\text{spec}}
\define\IA{\text{IA}}
\define\AIA{\text{AIA}}
\define\endpf{\hbox{\vrule height1.5ex width.5em}}
\define\Remark{{{\flushpar {\bf Remark.}\quad}}}

\define\alg{{{\frak g}}}
\define\dalg#1{{{\alg}^{(#1)}}}
\define\qalg{{\bar{\alg}}}
\define\center{{\frak z}}
\define\gp{{G}}
\define\qgp{{\bar{G}}}
\define\dgp#1{{\gp^{(#1)}}}
\define\qmetric{{\bar{g}}}

\define\GG{{\Gamma}}
\define\G#1{{{{\GG}_{#1}}}}
\define\qGG{{\bar{\GG}}}
\define\qG#1{{\bar\GG_{#1}}}

\define\Nilmfld{{\GG \backslash \gp}}
\define\nilmfld#1{{\G{#1} \backslash \gp}}

\define\qnilmfld#1{{\qG{#1} \backslash \bar{\gp}}}

\define\hh{{{\frak h}}}
\define\nn{{{\frak n}}}
\define\vv{{{\frak v}}}

\define\BX{{\bar{X}}}
\define\BY{{\bar{Y}}}
\define\BZ{{\bar{Z}}}
\define\BE{{\bar{E}}}
\define\BU{{\bar{U}}}
\define\BV{{\bar{V}}}
\redefine\b#1#2{{\bar{{#1}}_{#2}}}

\define\A#1#2#3{{A^{#1}_{{#2}{#3}}}}
\redefine\B#1#2#3{{B^{#1}_{{#2}{#3}}}}
\define\C#1#2#3{{C^{#1}_{{#2}{#3}}}}
\define\p#1#2{{\frac{\partial}{\partial {#1}_{#2}}}}
\redefine\d#1#2{{\dot{#1}_{#2}(s)}}

\topmatter

\title
The Marked Length Spectrum Versus the Laplace Spectrum on Forms on Riemannian Nilmanifolds
\endtitle

\rightheadtext{The Marked Length Spectrum on Nilmanifolds}

\author
Ruth Gornet
\endauthor

\affil
Texas Tech University
\endaffil

\address Ruth Gornet: Texas Tech University; 
Department of Mathematics; Lubbock, Texas \ 79409-1042;
May 1995
\endaddress

\keywords
Marked length spectrum, 
Length spectrum, 
Laplace spectrum, 
Laplace spectrum on forms,
Closed geodesics, 
Nilpotent Lie groups
\flushpar\indent Research at MSRI supported in part by NSF
grant DMS-9022140. Research at MSRI and Texas Tech supported in part
by NSF grant DMS-9409209
\endkeywords



\abstract
The subject of this paper is the relationship among
the marked length spectrum, the length spectrum,
the Laplace spectrum on functions,
and the Laplace spectrum on forms on Riemannian nilmanifolds.
In particular, we show that for a large class of three-step
nilmanifolds, if a pair of nilmanifolds in this class has the 
same marked length spectrum, they necessarily share the same
Laplace spectrum on functions. In contrast, we present the first
example of a pair of isospectral
Riemannian manifolds with the same
marked length spectrum but not the same spectrum on
one-forms. 
Outside of the standard spheres vs. the  Zoll spheres,
which are not even isospectral,
this is the only example of a pair of Riemannian manifolds
with the same marked length spectrum, but not the
same spectrum on forms.
This partially extends and partially contrasts
the work of Eberlein,
who showed that on two-step nilmanifolds, the same marked
length spectrum implies the same Laplace spectrum both on functions
and on forms.    
\endabstract
\endtopmatter


\document

\subheading{Section 1:  Introduction}
\bigskip

The {\it spectrum} of a closed Riemannian manifold $(M,g)$, 
denoted $\spec(M,g)$, is the collection of eigenvalues with multiplicities of the 
associated La\-place--Bel\-tra\-mi operator acting on smooth functions.
Two Riemannian manifolds $(M,g)$ and $(M',g')$ are said to be
{\it isospectral} if $\spec(M,g)=\spec(M',g').$

The Laplace--Beltrami operator may be extended to act on 
smooth $p$-forms by  $\Delta = d \delta + \delta d,$
where $\delta$ is the adjoint of $d$ and $p$ is 
a positive integer. We call its eigenvalue spectrum the
{\it p-form spectrum.}

The {\it length spectrum \/} of a Riemannian manifold 
is the set of lengths of smoothly closed geodesics, 
counted with multiplicity.  The multiplicity of a length 
is defined as the number of distinct 
free homotopy classes that contain 
a closed geodesic of that length. 
We denote the length spectrum of $(M,g)$ by 
$[L]\hy\spec(M,g).$
This is a natural notion, since the geodesic of shortest
length in a free homotopy class is just the shortest
loop representing that class.
(Note that other definitions of multiplicity appear in 
the literature.)

Two Riemannian manifolds $(M_1 ,g_1)$ and $(M_2,g_2)$
have the same {\it marked length spectrum \/} 
if there exists an isomorphism between the fundamental groups 
of $M_1$ and $M_2$
such that corresponding free homotopy classes
contain smoothly closed geodesics of the same length.
Clearly, manifolds with the same marked 
length spectrum  necessarily have the same length spectrum.

The purpose of this paper is to study the relationship among the
marked length spectrum, the length spectrum, the Laplace spectrum
on functions and the Laplace spectrum on forms on Riemannian
nilmanifolds. 

The relationship between the Laplace spectrum and lengths of 
closed geodesics arises from the study of the wave equation
(see \cite{DGu}, \cite{GuU}), and in the case of compact, hyperbolic
manifolds, from the Selberg Trace Formula (see \cite{C}, Chapter XI).
Colin de Verdiere \cite{CdV} has shown that generically,
the Laplace spectrum determines the length spectrum.
On Riemann surfaces, Huber showed that 
the length spectrum and the Laplace spectrum
are equivalent notions (see \cite{Bu} for an exposition). 

The Poisson formula gives the relationship between the
Laplace spectrum and length spectrum of flat tori,
with the result that
pairs of flat tori are isospectral if and only if
they share the same length spectrum (see \cite{CS}, \cite{G3}).
Pesce \cite{P2} has computed a Poisson-type formula relating 
the Laplace spectrum and length spectrum of Heisenberg manifolds,
and
has also shown that pairs of Heisenberg manifolds
that are isospectral must have the same lengths of closed geodesics. 
Previously, Gordon \cite{G1} exhibited the first examples of isospectral
manifolds that do not have the same length spectrum. 
These Heisenberg manifolds have the same lengths of closed geodesics.
However, the length spectra often differ in 
the multiplicities that occur. All known examples of
manifolds that are isospectral have the same lengths of closed geodesics. 

The marked length spectrum contains significantly more
geometric information than the length spectrum.
Croke \cite{Cr} and Otal \cite{Ot1}, \cite{Ot2} 
independently showed that if
a pair of compact
surfaces with negative curvature have the same marked
length spectrum, they are necessarily isometric. 
The same is true for flat tori (see \cite{G3}). 
Recently Eberlein \cite{E1} showed that for two-step nilmanifolds,
the same marked length spectrum implies the same Laplace spectrum
both on functions and on $p$-forms for all $p.$  (See Section 3
for more details.) 

However, the standard sphere and the Zoll sphere (see \cite{Bes})
have the same
marked length spectrum (trivially so, as they are both simply
connected and by definition have the same lengths of closed geodesics),
yet they are not even isospectral on functions.  
Indeed, any manifold isospectral to 
a standard sphere of dimension less than or equal to six must
be isometric to it (see \cite{B2}).
Examples of pairs of Riemannian manifolds that are isospectral
on functions but not on forms are sparse. 
Most constructions for producing pairs of isospectral manifolds
can be explained by Sunada's method \cite{S}  or its generalizations \cite{DG},
\cite{GW1}, \cite{B3}.
Pairs of manifolds constructed by the Sunada techniques necessarily
have the same $p$-form spectrum for all $p.$

For any choice of $P\in\Z^+,$ Ikeda \cite{I2} has constructed 
examples of isospectral lens spaces that are isospectral on $p$-forms for 
$p=0,1,\cdots,P$ but not isospectral on $(P+1)$-forms.
A straightforward argument shows that 
for the family of lens spaces considered by Ikeda, if a pair of
lens spaces
in this family 
has the same marked length spectrum, they are necessarily
isometric.
Gordon \cite {G2} has constructed pairs of Heisenberg manifolds
that are isospectral on functions, but not isospectral on one-forms.  
A consequence of Eberlein's theorem is that Heisenberg manifolds with
the same marked length spectrum
are necessarily isometric. (See Section 3 for more details.)
The only other known examples of manifolds that are isospectral
on functions but not isospectral on forms 
are pairs of isospectral three-step nilmanifolds
presented and studied in \cite{Gt3}.  These examples
are studied further here.

This paper focuses almost exclusively on three-step nilmanifolds.
The main results are a partial extension and a partial converse
to Eberlein's theorem for higher-step nilmanifolds.

\proclaim {Main Theorem 3.2.2}
For a large class of three-step nilmanifolds, if a pair
of nilmanifolds in this class has the same marked length spectrum,
they necessarily share the same Laplace spectrum on functions.
\endproclaim

\proclaim{Main Example}
Example V in the table below exhibits the first example
of
a pair of isospectral Riemannian manifolds with the same marked length
spectrum, but not the same spectrum on one-forms.
\endproclaim

These results have led to the following.

\proclaim{Conjecture}
Pairs of Riemannian nilmanifolds with the same marked
length spectrum are necessarily isospectral on functions.
\endproclaim

Background ideas and notation are established and explained
in Section 2.
In \cite{Gt3}, we presented a new construction for producing
pairs of isospectral nilmanifolds of arbitrary-step.
In Section 3, this construction together with Eberlein's 
theorem and techniques
from Riemannian geometry are used to prove the 
Main Theorem. 

Also in \cite{Gt3}, we presented new examples of isospectral
three-step nilmanifolds with combinations of
properties described in the table below.
For consistency, the numbering of the examples in this paper
coincides with the numbering of the examples in \cite{Gt3}.
Note that Example V is also the Main Example.

The spectrum on functions, spectrum on forms, quasi-regular
representations, and fundamental groups of these examples
were examined in \cite{Gt3}. In Section 4 and 5
we compare the length spectrum
and marked length spectrum of these examples. 
The pairs of isospectral 
manifolds described below have the same lengths of closed geodesics. 
However, the length spectra often differ in the multiplicities that occur.
 
\pagebreak

\centerline{Table I:  New Examples of Isospectral Manifolds}
\medskip
{\eightpoint
\centerline{$
\vbox{\offinterlineskip
\halign{\strut\vrule#&\ #\hfil\ &&\vrule#&\ \hfil#\hfil\ \cr
\noalign{\hrule}
&Pair of 3-Step&&$\forall p$ Same&&Rep. Equiv.&&Isomorphic&&Same&&Same&\cr
&Isospectral&&$p$-form&&Fundamental&&Fundamental&&Length&&Marked Length&\cr
&Nilmanifolds&&Spectrum&&Groups&&Groups&&Spectrum&&Spectrum&\cr
\noalign{\hrule}
&I(7 dim)&&Yes&&Yes&&No&&No&&No&\cr
\noalign{\hrule}
&II(5 dim)&&Yes&&Yes&&Yes&&Yes&&No&\cr
\noalign{\hrule}
&III$\backslash$IV(7$\backslash$5 dim)&&No&&No&&No&&No&&No&\cr
\noalign{\hrule}
&V(7 dim) &&No&&No&&Yes&&Yes&&Yes&\cr
\noalign{\hrule}}}
$}
}
\medskip

All of the examples described in the above table are of the form
$(\Nilmfld, g)$, where $\gp$ is a three-step nilpotent
Lie group, $\GG$ is a cocompact, discrete subgroup of $\gp$ (i.e.
$\Nilmfld$ compact) and $g$ arises from a left invariant
metric on $\gp.$
Two cocompact, discrete subgroups $\Gamma_1$ and 
$\Gamma_2$ of a Lie group $G$ are called {\it representation equivalent} if 
the associated quasi-regular representations are unitarily equivalent. 
If $\G1$ and $\G2$ are representation equivalent, then
$(\nilmfld1, g)$ and $(\nilmfld2, g)$ are necessarily isospectral
on functions and on smooth $p$-forms for any choice of left invariant metric
$g$ on $\gp.$

\Remark
Example I provided the first example of a pair of 
representation equivalent subgroups of a solvable Lie group producing
nilmanifolds with unequal length spectra.  This cannot
happen in the two-step nilpotent case.   
The relationship between the quasi-regular representation
and the length spectrum and marked length spectrum of 
nilmanifolds is studied in \cite{Gt2}, where we 
also present the first examples
of pairs of representation equivalent subgroups of
two-step nilpotent Lie groups that do not produce
nilmanifolds with the same marked length spectrum.
Example I is also the first example of a pair of
nonisomorphic, representation equivalent subgroups
of a solvable Lie group. See \cite{Gt1} for more details.
Note that nilpotent Lie groups are necessarily solvable.

Some of the contents of this paper are 
contained in the author's thesis
at Washington University in St. Louis in partial fulfillment of the requirements
for the degree of Doctor of Philosophy.  
The author wishes to express
deep gratitude to her advisor, Carolyn S. Gordon, for all of her suggestions,
encouragement, and support.
The author also wishes to thank Patrick Eberlein for helpful conversations.

\bigskip

\subheading{Section 2:  Background and Notation}
\bigskip

\subheading{Section 2.1: Definitions}
\bigskip

Let $\gp$ be a simply connected Lie group with Lie
algebra $\alg.$  
A metric on $\gp$ is {\it left invariant \/} if left translations 
are isometries.
Note that a left invariant metric is determined by 
a choice of orthonormal basis of the Lie algebra
$\alg$ of $\gp.$

Let $\GG$ be a a cocompact, 
discrete subgroup of $\gp.$
A left invariant metric 
$g$ on $\gp$ descends to a Riemannian metric on 
$\Nilmfld,$ which we also denote by $g.$ 
This paper focuses exclusively on manifolds of the form
$(\Nilmfld, g),$
where $g$ arises from a left invariant metric on $\gp.$
 
On manifolds of the form $(\Nilmfld,g),$
the Laplace--Beltrami operator is
$$\Delta = -\sum_{i=1}^n {E_i}^2,$$ where $\{E_1, \cdots, E_n\}$ 
is an orthonormal basis of the Lie algebra $\alg$ of $G.$

Recall that the free homotopy classes of a manifold $\Nilmfld$ 
correspond to the conjugacy classes in $\GG.$  
We will denote by $[\gamma]_\GG$ the free homotopy 
class of $\Nilmfld$ represented by $\gamma \in \GG.$ 
That is, 
$[\gamma]_\GG = \{ \hat{\gamma} \gamma \hat{\gamma}^{-1} : 
\hat{\gamma} \in \GG \}.$

We write $\lambda \in [\gamma]_\GG$ 
if there exists a closed geodesic of length 
$\lambda>0$ in the free homotopy class 
$[\gamma]_\GG$ of $(\Nilmfld, g).$

Let $\gamma$ be an element of $\GG.$ 
We say a geodesic $\sigma$ of $(\gp, g)$ 
is {\it translated by the element $\gamma$
with period $\lambda>0$}  if
$$\gamma \sigma(s) = \sigma(s + \lambda) \qquad \forall s \in \R.$$
If $\sigma$ is a unit speed geodesic,
then $\sigma$ projects to a closed geodesic on $(\Nilmfld, g)$ 
of length $\lambda,$  and $\sigma$ is contained
in the free homotopy class $[\gamma]_\GG.$

As the projection $(G,g) \rightarrow (\Nilmfld, g)$ 
is a Riemannian covering,
all closed geodesics of $(\Nilmfld, g)$ 
must arise in this fashion.
So to study the closed geodesics of $(\Nilmfld, g),$ 
it is enough to study the  $\gamma$-translated geodesics of $(G,g).$

Let $\sigma(s)$ be a geodesic of $\gp$ through 
$p = \sigma(0).$  Let $\hat{\sigma}(s) = p^{-1}\sigma(s).$
As left translations are isometries, 
$\hat{\sigma}$ is a geodesic of $\gp$ through $e.$
If $\sigma$ is translated by $\gamma$ with period $\lambda,$
then $\hat\sigma$ is translated by  $(p^{-1}\gamma p)$, also with
period $\lambda.$  
To see this, 
note that if $\gamma \sigma(s) = \sigma(s+\lambda),$  then 
$$(p^{-1} \gamma p)\hat{\sigma}(s) = 
(p^{-1}\gamma p)p^{-1} \sigma(s) = 
p^{-1}\gamma \sigma(s)  = 
p^{-1}\sigma(s+\lambda) = 
\hat{\sigma}(s+\lambda).$$

\flushpar {\bf 2.1.1 Notation.} In summary,  $\lambda \in [\gamma]_\Gamma$
if and only if there exists $x=p^{-1}\gamma p \in [\gamma]_G$ and a unit
speed geodesic $\sigma(s)$ on $(G,g)$ 
through $e=\sigma(0)$ such that 
$x \sigma(s) = \allowbreak \sigma(s+\lambda), \ \forall s \in \R.$  
That is, $x$ translates $\sigma$ with period $\lambda.$ 
Here $[\gamma]_G$ denotes the conjugacy class of $\gamma$
in $G.$

With this notation, a pair of manifolds 
$(\G1\backslash\gp_1, g_1)$ and 
$(\G2\backslash\gp_2, g_2)$ share the same marked length spectrum
if and only if there exists an isomorphism $\Phi:\G1 \rightarrow\G2$
such that for all $\gamma\in\G1,$
$$\lambda\in[\gamma]_\G1\text{ if and only if }\lambda\in[\Phi(\gamma)]_\G2.$$
We say that the isomorphism $\Phi$ {\it marks} the length spectrum 
between $(\G1\backslash\gp_1, g_1)$ and 
$(\G2\backslash\gp_2, g_2).$
\bigskip

\subheading{Section 2.2 Nilmanifolds}
\bigskip

Let $\alg$ be a Lie algebra.
We denote by $\dalg1$  
the derived algebra $[\alg,\alg]$ of $\alg.$ 
That is, $\dalg1$ is the Lie subalgebra of $\alg$ generated
by all elements of the form $[X,Y]$ for $X,Y$ in $\alg.$
Inductively, define $\dalg{k+1}=[\alg,\dalg{k}].$
The Lie algebra
$\alg$ is said to be {\it k-step nilpotent} if $\dalg{k} \equiv 0$
but $\dalg{k-1} \not \equiv 0.$
A Lie group $\gp$ is called {\it $k$-step nilpotent} if its Lie algebra is.

If $\gp$ is a nilpotent Lie group with cocompact, discrete subgroup $\GG,$
the locally homogeneous space $\Nilmfld$ is called a {\it nilmanifold.}
If $\gp$ is an abelian Lie group,
then $\GG$ is merely a lattice of rank $n$ in $G,$
where $n$ is the dimension of $G.$
In this case, $\log\GG$ is also a lattice in $\alg.$

Let $\exp$  denote  the Lie algebra 
exponential from $\alg$ to $G.$
The Campbell-Baker-Haus\-dorff formula gives us the group operation of 
$\gp$ in terms of $\alg.$  Namely, for $X,Y \in \alg:$
$$\exp(X)\exp(Y) = \exp(X+Y+\frac1{2}[X,Y]+\frac1{12}[X,[X,Y]]+
\frac1{12}[Y,[Y,X]]+\cdots),$$
where the remaining terms are higher-order brackets.
Note that for two-step nilpotent Lie groups, 
only the first three terms in the
right-hand side are nonzero.  For three-step groups, only the first
five terms are nonzero.
If $\alg$  is nilpotent and $\gp$ is simply connected,
then $\exp$ is a diffeomorphism from $\alg$
onto $\gp.$  Denote its inverse by $\log.$

If $\gp_1$ and $\gp_2$ are nilpotent Lie groups with cocompact,
discrete subgroups $\G1$ and $\G2,$ respectively, any
abstract group isomorphism 
$\Phi:\G1\rightarrow\G2$ lifts uniquely to a Lie group 
automorphism $\Phi:\gp_1\rightarrow\gp_2.$

For details of cocompact, discrete sugbroups of nilpotent Lie 
groups, see \cite{Ra}.

\proclaim{2.2.1 Definition}
Let $\Phi$ be a Lie group automorphism of $G.$
Let $\GG$ be a cocompact, discrete subgroup of $G.$

(i) We call $\Phi$ an  {\it almost inner automorphism} if 
for all elements $x$ of $\gp$
there exists $a_x$ in $\gp$ such that   $\Phi(x)=a_xxa_x^{-1}.$

(ii) We say $\Phi$ is a {\it $\GG$-almost inner automorphism} if 
for all elements $\gamma$ of $\GG$
there exists $a_\gamma$ in $\gp$ such that
$\Phi(\gamma)=a_\gamma \gamma a_\gamma^{-1}.$
\endproclaim

Denote by $\IA(\gp)$ (respectively, $\AIA(\gp),\GG\hy\AIA(\gp)$ )
the group of inner automorphisms (respectively,
almost inner automorphisms, 
$\GG$-almost inner automorphisms)
of $\gp.$
Note that $\IA(\gp) \subset \AIA(\gp) \subset \GG\hy\AIA(\gp).$

\proclaim{2.2.2 Theorem (Gordon and Wilson, Gordon \cite{GW1},\cite{G1}) }
Let $\gp$ be an exponential solvable Lie group,
and let $\G1$ and $\G2$ be cocompact, discrete subgroups
of $G.$  Let $\Phi$ be a $\G1$-almost inner
automorphism of $\gp$ such that $\Phi(\G1)=\G2.$
Then $(\nilmfld1,g)$ and $(\nilmfld2,g)$
are isospectral on functions and on forms
for any choice of left invariant metric $g$ on $\gp.$
Moreover,  
the automorphism $\Phi$ marks the length spectrum 
between $(\nilmfld1,g)$ and $(\nilmfld2,g).$
\endproclaim

Note that a nilpotent Lie group is 
necessarily exponential solvable.
\bigskip

\subheading{Section 3: The Marked Length Spectrum vs. 
the Laplace Spectrum on Functions of Three-Step Nilmanifolds }
\bigskip

Throughout this section, $\gp$ is a simply connected, 
$k$-step nilpotent Lie group, with Lie algebra $\alg,$  
$\GG$ is a cocompact, discrete subgroup of $\gp,$ and $g$ is a 
left invariant metric on $\gp$ which descends to a metric on
$\Nilmfld,$ also denoted by $g.$ 
We denote the center of $\alg$ by $\center$
and the center of $\gp$ by $Z(\gp).$
Let $L_x$ denote 
left multiplication by $x\in\gp.$  As $g$ is left invariant, 
$L_x$ is always an isometry of $(\gp,g).$
Let $\dgp{k} = \exp(\dalg{k})$ denote the $k$th derived subgroup of $\gp.$ 
Note that if $\gp$ is $k$-step nilpotent, then $\dgp{k-1} \subset Z(\gp).$
\bigskip

\subheading{\S3.1 Preliminaries}
\bigskip

\proclaim{3.1.1 Theorem}
Let $\gp$ be a three-step nilpotent Lie group with left invariant
metric $g.$ 
Let $\sigma$ be a geodesic on $(\gp,g)$ 
that is translated by the element $\gamma \in \gp$
with period $\lambda > 0.$  Let $p = \sigma(0).$  
Then 
$$\left<L_{p*}\left(\left[\log\left(p^{-1}\gamma p\right),\alg\right]\right),
\dot{\sigma}(0)\right>_p \equiv 0.$$
\endproclaim
\medskip

\Remark
This is the three-step generalization 
of a result due to Eberlein \cite{E1}.
Recently Dorothee Schueth \cite{Sch} has given an 
elegant proof, which generalizes the result to nilpotent Lie 
groups of {\it arbitrary} step.
\medskip

\demo{Outline of Proof of 3.1.1}

We briefly describe the basic steps in the original three-step proof.
For details, see \cite{Gt4}, Chapter 4.

Let $\gp$ be a simply connected, three-step nilpotent Lie group with Lie algebra $\alg$ and left invariant metric $g.$
Let $\alg = \nu \oplus \dalg1,$ 
where $\nu$ is the orthogonal complement of $\dalg1$ in $\alg.$  
Let $\dalg1 = \zeta \oplus \dalg2,$ 
where $\zeta$ is the orthogonal complement of 
$\dalg2$ in $\dalg1.$   
Thus $\alg = \nu \oplus \zeta \oplus \dalg2.$

Let $\{ X_1, X_2, \cdots, X_J\}$ be an orthonormal basis of $\nu.$
Let $\{ Z_1, Z_2, \cdots, Z_K\}$ be an orthonormal basis of $\zeta,$
and let $\{ W_1, W_2, \cdots, W_T\}$ be an orthonormal basis of $\dalg2.$
Throughout this proof
the indices $i,j,$ and $l$ run from $1$ to $J,$ 
the indices $h$ and $k$ run from $1$ to $K,$ 
and the indices $t$ and $r$ run from $1$ to $T.$

Define $\A{k}{i}{j}, \B{t}{i}{j}, \C{t}{i}{k}$ by
$$\align [X_i, X_j] &= \sum_k\A{k}{i}{j}Z_k + 
\sum_t\B{t}{i}{j}W_t\\
[X_i, Z_k] &= - [Z_k, X_i] =  \sum_t\C{t}{i}{k}W_t.
\endalign$$
As $[X_i,X_j] = - [X_j,X_i],$ 
we have $\A{k}{i}{j}=-\A{k}{j}{i}$ 
and $\B{t}{i}{j}=-\B{t}{j}{i}.$ 
By the Jacobi equation $[\dalg1, \dalg1] \subset [\alg, \dalg2]\equiv 0.$ 
Thus $[Z_k, Z_h] = 0.$
Finally, by applying the Jacobi equation to $X_i, X_j, X_k$
and examining the $W_t$ coefficient,
we obtain:
$$0 = \sum_k\left(\A{k}{j}{l}\C{t}{i}{k} + \A{k}{i}{j}\C{t}{l}{k} +
\A{k}{l}{i}\C{t}{j}{k}\right).$$

For Lie algebras with a left invariant metric, 
the covariant derivatives can be calculated via

\centerline{$<\nabla_VY,U> =\frac1{2} <[U,V],Y> + 
\frac1{2}<[U,Y],V> + 
\frac1{2}<[V,Y],U>$}

\flushpar for $U,V,Y$ in $\alg.$
We obtain the covariant derivatives:
$$\split 
\nabla_{X_i}X_j &= \frac1{2}\sum_{k}\A{k}{i}{j}Z_k + \frac1{2}\sum_{t}\B{t}{i}{j}W_t, \\
\nabla_{X_i}Z_k &= \frac1{2}\sum_{j}\A{k}{j}{i}X_j + \frac1{2}\sum_{t}\C{t}{i}{k}W_t, \\
\nabla_{Z_k}X_i &= \frac1{2}\sum_{j}\A{k}{j}{i}X_j - \frac1{2}\sum_{t}\C{t}{i}{k}W_t, \\
\nabla_{X_i}W_t &= \nabla_{W_t}X_i = \frac1{2}\sum_{j}\B{t}{j}{i}X_j - \frac1{2}\sum_{k}\C{t}{i}{k}Z_k, \\
\nabla_{Z_k}{Z_h} &=\nabla_{W_t}{W_r} = 0, \\
\nabla_{Z_k}{W_t} &=\nabla_{W_t}{Z_k} = \frac1{2}\sum_{j}\C{t}{j}{k}X_j.
\endsplit$$

For $x\in\gp,$   
$x=\exp\left(\sum_j x_jX_j+\sum_k z_kZ_k+\sum_t w_tW_t\right)$
gives us a global coordinate system on $\gp.$
With this coordinate system, a straightforward computation
shows us that
$$\align
X_j &= \p{x}{j} + 
\sum_{k}\left(\frac1{2}\sum_{i}x_i\A{k}{i}{j}\right)\p{z}{k}\\
& \qquad + 
\sum_{t}\left(\frac1{2}\sum_{i}x_i\B{t}{i}{j} - 
\frac1{2}\sum_{k}\C{t}{j}{k}z_k + \frac1{12}\sum_{i,l,k}x_i\C{t}{i}{k}x_l\A{k}{l}{j}\right)\p{w}{t},\\
Z_k& = \p{z}{k} + \sum_t\left(\frac1{2}\sum_{i}x_i\C{t}{i}{k}\right) \p{w}{t},\\
W_t &= \p{w}{t}.
\endalign$$

Let 
$\sigma(s) = 
\exp\left(\sum_j{x_j(s)} X_j + 
\sum_k z_k(s) Z_k + 
\sum_t w_t(s) W_t\right)$ be a geodesic
of $(\gp, g)$ with
initial velocity 
$\dot\sigma(0)=\sum_j\bar{x}_j X_j + 
\sum_k \bar{z}_k Z_k + 
\sum_t \bar{w}_t W_t.$ 
A straightforward computation 
of $\nabla_{\dot\sigma(s)}\dot\sigma(s)\equiv0$ produces
the following geodesic equations for a 
three-step nilpotent Lie group,
reduced to a system of $n$-ordinary
differential equations.

$$\align
\d{x}{j} &= 
-\sum_{l,k}x_l(s)\A{k}{j}{l}\b{z}{k} - 
\sum_{l,t}x_l(s)\B{t}{j}{l}\b{w}{t} - 
\sum_{k,t}z_k(s)\C{t}{j}{k}\b{w}{t}\\
&\qquad  -
\frac1{2}\sum_{i,l,k,t}x_i(s)x_l(s)\b{w}{t}\C{t}{i}{k}\A{k}{j}{l} + 
\b{x}{j}\\
\\
\d{z}{k} &= 
\frac1{2}\sum_{i,j}x_i(s)\d{x}{j}\A{k}{i}{j} + \sum_{j,t}x_j(s)\b{w}{t}\C{t}{j}{k} + 
\b{z}{k}\\
\\
\d{w}{t} &= 
\frac1{2}\sum_{i,j}x_i(s)\d{x}{j}\B{t}{i}{j} - \frac1{2}\sum_{j,k}\d{x}{j}z_k(s)\C{t}{j}{k} + \frac1{2}\sum_{j,k}x_j(s)\d{z}{k}\C{t}{j}{k}\\
& \qquad -
\frac1{6}\sum_{i,j,k,l}x_i(s)\d{x}{j}x_l(s)\C{t}{i}{k}\A{k}{l}{j} + 
\b{w}{t}
\endalign$$

If we assume that a geodesic $\sigma(s)$ starts at the identity
and is translated by the element $\gamma,$
then a lengthy but
straightforward (brute-force) calculation yields
$$\left<\left[\log(\gamma),\alg\right], 
\dot\sigma(0)\right>_e \equiv 0.$$
Here one uses the extensively the fact that
if $\gamma \sigma(s)=\sigma(s+\lambda),$
then $L_{\gamma*}\left(\dot\sigma(s)\right)=\dot\sigma(s+\lambda).$

In the general case, let $\sigma(s)$ be a geodesic of $\gp$ through 
$p = \sigma(0).$  Let $\alpha(s) = p^{-1}\sigma(s).$
Then $\alpha$ is a geodesic of $\gp$ through $e.$
If $\sigma$ is translated by $\gamma$ with period $\lambda,$
then $\alpha$ is translated by $p^{-1}\gamma p,$ also
with period $\lambda.$
Thus
$$\left<\left[\log\left(p^{-1} \gamma p\right), \alg \right], 
\dot\alpha(0) \right>_e \equiv 0.$$

But $\dot\alpha(0) = (L_{p^{-1}})_*\left(\dot\sigma(0)\right).$
As our metric is left invariant, we obtain
$$\left<L_{p*}\left(\left[\log\left(p^{-1} \gamma p\right), \alg\right]\right),
\dot{\sigma}(0)\right>_p = 0,$$
as desired.\endpf
\enddemo
\medskip

\Remark Ron Karidi \cite{K} has recently given a formulation of 
the geodesic equations for an arbitrary nilpotent Lie group with
a left invariant metric.  As above, this formulation is 
in terms of an orthonormal basis and structure constants of 
the Lie algebra.
\medskip

\flushpar {\bf 3.1.2 Notation.} 
Let $\pi$ denote the projection from $\gp$ onto 
$\qgp = \gp/\dgp{k-1}.$ 
For $\GG$ a cocompact, discrete subgroup of $\gp,$
denote by $\qGG$ the image of $\GG$ under the canonical projection 
from $\gp$ onto $\qgp.$
The group $\qGG$ is then a cocompact, discrete subgroup of $\qgp.$
Let $\qmetric$ denote the metric on $\qgp$ defined by 
restricting the left invariant metric $g$ to an orthogonal 
complement of $\dalg{k-1}\subset\center,$ where $\alg$ is the 
Lie algebra of $\gp.$
With this choice of metric $\qmetric$ on $\qgp,$  
the mapping 
$$\pi:(\gp, g) \rightarrow(\qgp,\qmetric)$$ is a 
Riemannian submersion with totally geodesic fibers.

If $\Phi:\gp_1\rightarrow\gp_2$ is a Lie group mapping,
then necessarily $\Phi:\gp^{(k-1)}_1\rightarrow\gp^{(k-1)}_2.$
Let $\bar\Phi$ denote the canonical projection of $\Phi$
onto $\bar\Phi=\pi\circ\Phi:\qgp_1\rightarrow\qgp_2.$

The Lie algebra of $\qgp$ is $\qalg = \alg/\dalg{k-1}.$
We denote elements of 
$\qalg$ by $\BU$ where 
$\BU$ is the image of $U\in\alg$ 
under the canonical projection from 
$\alg$ onto $\qgp .$  Similarly, 
we will denote elements of $\qgp$ 
by $\bar{x}$ where $\bar{x}$ is the image of 
$x\in \gp$ under the canonical projection from 
$\gp$ onto $\qgp.$

All of the nilpotent Lie groups studied here have the 
following property.

\proclaim{3.1.3 Definition} Let $\gp$ be a simply connected, 
$k$-step nilpotent Lie group.  We say $\gp$ is 
{\it strictly nonsingular} 
if the following property holds:  
for all $z$ in $Z(\gp)$ and for all 
noncentral $x$ in $\gp$ there exists $a$ in $\gp$ such that 
$$[a,x] = z.$$
Here $[a,x]=axa^{-1}x^{-1}.$
Equivalently, the Lie algebra $\alg$ is 
{\it strictly nonsingular} if for all noncentral $X$ in $\alg,$  
$$\center \subset ad(X)(\alg).$$
That is, for all $X$ in $\alg-\center$ and all $Z$ in 
$\center$ there exists $Y$ in $\alg$ such that $[X,Y]=Z.$
\endproclaim

Note that for strictly nonsingular nilpotent Lie algebras,
$\center = \dalg{k-1}.$

\proclaim{3.1.4 Corollary}
Let $\gp$ be a simply connected, strictly nonsingular 
three-step nilpotent Lie group with left invariant
metric $g.$    Consider the Riemannian submersion
$(\gp,g)\rightarrow(\qgp,\qmetric).$ If 
$\sigma$ is a geodesic on $\gp$ such that 
$\gamma \sigma(s) = \sigma(s + \lambda)$ for some noncentral 
$\gamma$ in $\gp$ and some $\lambda > 0,$  then 
$\sigma$ is a horizontal geodesic.  That is, 
$$\left<{L_{\sigma(s)*}}(\center),\dot{\sigma}(s)\right> \equiv 0 
\qquad \qquad \forall s \in \R.$$  
\endproclaim

Before proving Corollary 3.1.4, recall the
following properties of Riemannian submersions.

\proclaim{3.1.5 Proposition (see \cite{GHL})}
Let $(M,g)\rightarrow(\bar{M},\qmetric)$ be a Riemannian
submersion.

(i)  Let $\alpha$ be a geodesic of $(M,g).$  
If the vector 
$\dot{\alpha}(0)$ is horizontal, then $\dot{\alpha}(s)$ 
is horizontal for all $s,$ and the curve $\pi \circ \alpha$ 
is a geodesic of $(\bar{M}, \qmetric)$ of the same length as $\sigma.$

(ii)  Conversely, let $p\in M$ and let $\sigma$ be a geodesic of 
$(\bar{M}, \qmetric)$ with $\sigma(0) = \pi(p).$  Then there exists a 
unique local horizontal lift $\hat{\sigma}$ of $\sigma$ through
$p=\hat\sigma(0),$  and 
$\hat{\sigma}$ is also a geodesic of $(M,g).$  
\endproclaim

\demo{Proof of Corollary 3.1.4}

By Theorem 3.1.1
$$\left<L_{p*}\left(\left[\log(p^{-1} \gamma p), \alg\right]\right),
\dot{\sigma}(0)\right>_p \equiv 0,$$
where $p = \sigma(0).$
By strict nonsingularity
$$\center = \dalg2 \subset \left[\log(p^{-1} \gamma p), \alg\right].$$
Thus $$\left<L_{p*}(\center),\dot{\sigma}(0)\right>_p \equiv 0.$$
Thus $\dot{\sigma}(0)$ is horizontal.
By Proposition 3.1.5, we know that 
$\dot{\sigma}(s)$ is horizontal for all $s \in \R.$
$\endpf$
\enddemo
\bigskip

\subheading{\S3.2 Main Theorem}
\bigskip

On two-step nilmanifolds, we have the following relationship betweeen the 
marked length spectrum and the $p$-form spectrum.

\proclaim{3.2.1 Theorem (Eberlein \cite {E})}
Let $\G1,\G2$ be cocompact, discrete subgroups of simply connected,
two-step nilpotent Lie groups $\gp_1,\gp_2$ with left invariant metrics 
$g_1 , g_2$ respectively.
Assume that  $(\G1\backslash\gp_1,g_1)$ and 
$(\G2\backslash \gp_2,g_2)$ 
have the same marked length spectrum, and
let $\Phi:\G1\rightarrow\G2$ be an isomorphism inducing
this marking.
Then $\Phi = {(\Phi_1\circ\Phi_2)\vert}_\G1,$ 
where $\Phi_2$ is a $\G1$-almost-inner automorphism of $\gp_1$,
and $\Phi_1$ is an isomorphism of $(\gp_1,g_1)$ onto $(\gp_2,g_2)$
that is also an isometry. Moreover, this factorization is unique.
In particular, $(\G1\backslash\gp_1,g_1)$ and
$(\G2\backslash\gp_2,g_2)$
have the same spectrum of the Laplacian on functions
and on $p$-forms for all $p.$
\endproclaim

\Remark Note that if $\GG\hy\AIA(\gp)=\IA(\gp),$
then the elements of $\GG\hy\AIA(\gp)$ are isometries
of $(\gp,g),$  where $g$ is any choice of left invariant metric 
$g$ of $\gp.$  
So by Theorem 3.2.1, 
any two-step nilmanifold with the same marked length 
spectrum as $(\Nilmfld,g)$
is necessarily isometric to it.
Note that this property applies to Heisenberg groups.
Thus pairs of Heisenberg manifolds with the same marked
length spectrum are necessarily isometric.
\medskip

We may now state the main result of this paper.

\proclaim{3.2.2 Main Theorem}
Let $\gp$ be a simply connected, strictly nonsingular, 
three-step nilpotent Lie group.  Let $\G1$ and $\G2$ be cocompact, 
discrete
subgroups of $\gp$  such that 
$\G1 \cap Z(\gp) = \G2 \cap Z(\gp).$ 
If $(\nilmfld1, g)$ and $(\nilmfld2, g)$ 
have the same marked length spectrum, 
then $(\nilmfld1, g)$ and $(\nilmfld2, g)$ 
are isospectral on functions.
\endproclaim
\medskip

To prove Theorem 3.2.2, we need the following.

\proclaim{3.2.3 Theorem \cite{Gt3, Theorem 3.2}}
Let $\gp$ be a simply connected, 
strictly nonsingular nilpotent Lie group with left 
invariant metric $g.$  
If \ $\G1$ and $\G2$ are cocompact, discrete subgroups 
of $\gp$ such that 
$$\G1 \cap Z(\gp) = \G2 \cap Z(\gp) \quad \text{and} 
\quad \spec\left(\qnilmfld1,\qmetric\right) = 
\spec\left(\qnilmfld2,\qmetric\right),$$ then  
$$\spec\left(\nilmfld1,g\right) = \spec\left(\nilmfld2,g\right).$$
\endproclaim
\medskip

\proclaim{3.2.4 Theorem}
Let $\gp$ be a simply connected, strictly nonsingular 
three-step nilpotent Lie group with cocompact, discrete 
subgroup $\GG$ and left invariant metric $g.$  
Let $\gamma$ be a {\it noncentral} element of $\GG.$ 
Then we have the following condition: 
$$\lambda \in [\gamma]_\GG \quad  \text{if and only if }\quad 
\lambda \in [\pi(\gamma)]_\qGG.$$
\endproclaim
\medskip

Assume for the moment that Theorem 3.2.4 is true.

\proclaim{3.2.5 Corollary}  Let $\gp_1$ and $\gp_2$ be simply connected, 
strictly nonsingular, three-step nilpotent 
Lie groups with cocompact, discrete subgroups   
$\G1$ and $\G2$ and left invariant metrics 
$g_1$ and $g_2,$ respectively.  Let 
$\Phi$ mark the length spectrum between 
$(\G1 \backslash \gp_1, g_1)$ and 
$(\G2\backslash \gp_2, g_2).$  Then 
$\bar\Phi$ must mark the length spectrum between 
$(\qG1 \backslash \qgp_1, \qmetric_1)$ and 
$(\qG2 \backslash \qgp_2, \qmetric_2).$  
\endproclaim
\medskip

\demo{Proof of Corollary 3.2.5}

Let $\lambda \in [\pi(\gamma)]_\qG1, \pi(\gamma)\neq 0.$
By (3.2.4) $\lambda \in [\gamma]_\G1.$
By hypothesis $\lambda \in [\Phi(\gamma)]_\G2.$
By (3.2.4)  again
$\lambda \in [\pi(\Phi(\gamma))]_\qG2 = 
[\bar\Phi(\pi(\gamma))]_\qG2.$  

Reversing the roles of $\qG1$ and $\qG2,$ 
we obtain the desired result.
$\endpf$
\enddemo
\medskip

\demo{Proof of Main Theorem 3.2.2}

Let $\Phi$ mark the length spectrum between $(\nilmfld1, g)$ and 
$(\nilmfld2, g).$
By (3.2.5) we know that $\bar\Phi$ must mark the 
length spectrum between $\left(\qnilmfld1, \qmetric\right)$ and 
\newline
$\left(\qnilmfld2, \qmetric\right).$
By Theorem 3.2.1
$\spec(\qnilmfld1, \qmetric)=\spec(\qnilmfld2, \qmetric).$

The result now follows directly from Theorem 3.2.3.
$\endpf$ 
\enddemo
\medskip

It remains only to prove Theorem 3.2.4, which 
follows directly from the following two lemmas.

\proclaim{3.2.6 Lemma}
Let $\gp$ be a simply connected, strictly nonsingular 
three-step nilpotent Lie group with cocompact, 
discrete subgroup $\GG$ and left invariant metric $g.$  
Let $\gamma$ be a noncentral element of $\GG.$  
With the above notation, if $\lambda \in [\gamma]_\GG$ then 
$\lambda \in [\pi(\gamma)]_\qGG.$
\endproclaim

\demo{Proof of Lemma 3.2.6}

If $\lambda \in [\gamma]_\GG,$  then there exists a 
unit speed geodesic 
$\sigma(s)$ of $\gp$ through $e$   
such that $$p^{-1}\gamma p\sigma(s) = \sigma(s+\lambda)$$ 
for some $p \in\gp.$

By (3.1.4), $\sigma(s)$ is a horizontal geodesic,
and by (3.1.5), $\pi\circ\sigma(s)$ is a unit speed geodesic of 
$(\qgp, \qmetric).$ 

But 
$\pi(p^{-1}\gamma p\sigma(s)) = 
\pi(p^{-1})\pi(\gamma)\pi(p)\pi(\sigma(s)) = 
\pi(\sigma(s+\lambda)).$  
Thus $\pi(\sigma)$ is a unit speed geodesic translated by
$\pi(p^{-1})\pi(\gamma)\pi(p)$ with period $\lambda.$
That is, 
$\lambda \in [\pi(\gamma)]_\qGG, $ as desired. $\endpf$
\enddemo

\proclaim{3.2.7 Lemma}
For $\gp$ a simply connected, strictly nonsingular 
$k$-step nilpotent Lie group. Using the above notation, 
let $\lambda \in [\bar\gamma]_\qGG,$ 
where $\bar\gamma \neq e.$  
Then $\lambda \in [\gamma]_\GG$ for all 
$\gamma \in \pi^{-1}(\bar{\gamma}).$
\endproclaim
\medskip

\demo{Proof of Lemma 3.2.7}

Let $\sigma$ be a unit speed geodesic of 
$(\qgp, \qmetric)$  through $\bar{e} = \sigma(0)$ and 
translated by $p^{-1}\bar\gamma p$ for some $p \in \qgp.$

By (3.1.5), the unique horizontal lift 
$\hat{\sigma}$ of $\sigma$ with 
$\hat{\sigma}(0)=e$ is a geodesic of $(\gp,g).$

As both $\gp$ and $\qgp$ are complete, 
we see that $\hat{\sigma}$ is defined for all $s \in \R.$
We also have 
$\pi\circ\hat{\sigma}(s) = \sigma(s)$ for all $s \in \R.$  
To see this, note that the set $S$ of 
all such $s$ is nonempty as $0 \in S,$ open by completeness, 
and closed by uniqueness and smoothness.  Thus, $S=\R.$

Now $\pi(\hat{\sigma}(\lambda)) = p^{-1}\bar\gamma p.$   
Let $\hat{p}$ be such that $\pi(\hat{p}) = p.$  

Let $\gamma \in \pi^{-1}(\bar\gamma).$
Then 
$\pi(\hat{p}^{-1}\gamma\hat{p}) = 
p^{-1}\bar\gamma p = 
\pi(\hat{\sigma}(\lambda)).$  
Thus $(\hat{\sigma}(\lambda))(\hat{p}^{-1}\gamma\hat{p})^{-1}$ 
is a central element of $\gp.$

By strict nonsingularity, 
there exists $x \in \gp$ such that 
$$x^{-1}(\hat{p}^{-1}\gamma\hat{p})x(\hat{p}^{-1}\gamma\hat{p})^{-1} = 
\hat{\sigma}(\lambda)(\hat{p}^{-1}\gamma\hat{p})^{-1},$$  
that is $x^{-1}(\hat{p}^{-1}\gamma\hat{p})x= \hat{\sigma}(\lambda).$

If we let $p' = \hat{p}x,$  
then $\hat{\sigma}(\lambda)= {p'}^{-1}\gamma p'.$
Note that $\pi({p'}^{-1}\gamma p') = \pi(\hat\sigma(s)) = p^{-1}\bar\gamma p.$

We now show that 
${p'}^{-1}\gamma p'{\hat\sigma}(s) = 
\hat{\sigma}(s+\lambda) $ for all $s \in \R.$
Let  
$$\alpha(s) = ({p'}^{-1}\gamma p') ^{-1}\hat{\sigma}(s+\lambda).$$
Now $\alpha(0) = ({p'}^{-1}\gamma p')^{-1}\hat{\sigma}(\lambda) = e.$ 
Also, $\alpha(s)$ is horizontal since $g$ 
is left invariant and $\alpha$ is just a 
left translate of the horizontal curve 
$\hat{\sigma}.$ 
Moreover, 
$$\align
\pi(\alpha(s)) &= \pi(({p'}^{-1}\gamma p')^{-1}\hat\sigma(s+\lambda) ) = 
p^{-1}\bar\gamma^{-1}p\sigma(s + \lambda) \\
&= 
p^{-1}\bar\gamma^{-1} p p^{-1} \bar\gamma p \sigma(s) = 
\sigma(s).
\endalign$$
Thus $\alpha$ is a horizontal geodesic through 
$e \in G$ whose projection agrees with $\sigma.$  
By uniqueness in Proposition 3.1.5, 
$\alpha(s) = \hat\sigma(s) \quad \forall s \in \R.$ 

Consequently, 
$${{p'}^{-1}\gamma p'} \hat{\sigma}(s) = 
\hat{\sigma}(s+\lambda) $$ 
for all 
$s \in \R.$
Thus $$\lambda \in [\gamma]_\GG,$$ as desired.
\endpf
\enddemo
\bigskip

\subheading{\S3.3 Three-step Nilmanifolds with a One-Dimensional Center}

\proclaim{3.3.1 Theorem} Let $\gp$ be a simply connected, 
strictly nonsingular, three-step nilpotent Lie group 
with a one-dimensional center.  Let $\G1$ and $\G2$ 
be cocompact, discrete subgroups of $\gp$ such that 
$\G1 \cap Z(\gp) = \G2 \cap Z(\gp).$ Let $g$ be any 
left invariant metric on $\gp.$ Then 
$(\nilmfld1,g)$ and $(\nilmfld2,g)$ 
have the same marked length spectrum 
if and only if 
there exists an isomorphism 
$\Phi:\G1 \rightarrow \G2$ such that 
$\bar\Phi:\qG1 \rightarrow \qG2$ 
marks the length spectrum between 
$(\qnilmfld1, \qmetric)$ and $(\qnilmfld2, \qmetric).$ 
\endproclaim

\demo{Proof of Theorem 3.3.1}

The forward direction follows immediately from Corollary 3.2.5.

For the converse direction, 
assume that there exists an isomorphism 
$\Phi:\G1 \rightarrow \G2$ such that 
$\bar\Phi$ marks the length spectrum between 
$(\qnilmfld1, \qmetric)$ and $(\qnilmfld2, \qmetric).$

We need to show that for all $\gamma \in \G1,$ 
{$\lambda \in [\gamma]_\G1$ 
if and only if 
$\lambda \in [\Phi(\gamma)]_\G2.$}

We consider two cases:

Case 1:  $\gamma \in \G1 \cap Z(\gp).$

If $\lambda \in [\gamma]_\G1,$  
then there exists a geodesic 
$\sigma(s)$ of $\gp$ such that 
$\gamma \sigma(s) = \sigma(s+\lambda).$  

As $\Phi$ is an isomorphism, we know that 
$\Phi(\G1 \cap Z(\gp)) = \G2 \cap Z(\gp) = \G1 \cap Z(\gp),$  
and hence, $\Phi$ must map a generator of 
$\G1 \cap Z(\gp)$ into a generator of 
$\G1 \cap Z(\gp).$
There are only two such generators. 
Thus for all $\gamma \in \G1 \cap Z(\gp),$
either $\Phi(\gamma)=\gamma$ or $\Phi(\gamma)= \gamma^{-1}.$

Hence
$[\Phi(\gamma)]_\G2 = [\gamma]_\G2$ or 
$[\Phi(\gamma)]_\G2 = [\gamma^{-1}]_\G2.$

If $[\Phi(\gamma)]_\G2 = [\gamma]_\G2,$ 
then the geodesic $\sigma(s)$ of $\gp$ 
projects to a closed geodesic of $(\nilmfld2,g)$ 
of length $\lambda$ in the free homotopy class $[\gamma]_\G2.$

If $[\Phi(\gamma)]_\G2 = [\gamma^{-1}]_\G2,$ 
then the geodesic $\alpha(s) = \sigma(-s)$ of 
$\gp$ projects to a closed geodesic of 
$(\nilmfld2,g)$ of length $\lambda$ in $[\gamma^{-1}]_\G2.$

This argument also works for 
$\Phi^{-1}:\G2 \rightarrow \G1,$ 
which must necessarily mark the length spectrum.
Consequently, for all  $\gamma \in \G{1} \cap Z(\gp),$ 
$$\lambda \in [\gamma]_\G1 
\text{ if and only if } 
\lambda \in [\Phi(\gamma)]_\G2.$$

Case 2:  $\gamma \not \in Z(\gp)$

Let $\lambda \in [\gamma]_\G1.$  By strict nonsingularity and 
Theorem 3.2.4, 
we know that $\lambda \in [\pi(\gamma)]_\qG1.$
By assumption (ii), we know that 
$\lambda \in [\bar\Phi(\pi(\gamma))]_\qG2.$
Now $\pi(\Phi(\gamma)) = \bar\Phi(\pi(\gamma)).$  
Thus by Theorem 3.2.4 again we know 
$\lambda \in [\Phi(\gamma)]_\G2.$  
Reversing the roles of $\G1$ and $\G2$ 
in the above, we see that for all $\gamma \in \G1$ 
$\gamma\not\in\G1,$
$$\lambda \in [\gamma]_\G1 
\text{ if and only if } 
\lambda \in [\Phi(\gamma)]_\G2,$$
as desired. $\endpf$
\enddemo

\bigskip

\subheading{\S4 The Marked Length Spectrum vs. the One-Form  Spectrum}
\bigskip

The example below is the first example of a pair of isospectral
Riemannian
manifolds with the same marked length spectrum, but not the
same spectrum on one-forms. Outside of the standard vs. 
Zoll spheres, which are not even isospectral for dimension less
than or equal to six, this is the only example
of a pair of Riemannian manifolds that have the same marked length
spectrum but not the same spectrum on one-forms.
\medskip

\subheading{Example V}
\bigskip

We use the notation of Section 3. 

Consider the simply connected, strictly nonsingular, three-step 
nilpotent Lie group $\gp$ with Lie algebra
$$\alg = span_{\R}\{X_1, X_2, Y_1, Y_2, Z_1, Z_2, W \}$$
and Lie brackets
$$[X_1,Y_1] = [X_2, Y_2] = Z_1$$
$$[X_1,Y_2] = Z_2$$
$$[X_1,Z_1] = [X_2, Z_2] = [Y_1, Y_2] = W$$
and all other basis brackets zero.
\bigskip

We fix a left invariant metric on $\gp$ by letting 
$\{E_1,E_2,E_3,E_4,E_5,E_6,E_7\}$ be an orthonormal basis of $\alg$ where
$$
\align
E_1&= X_1-\frac{1}{2}X_2-\frac{1}{4}Y_2,\\
E_2&= X_2-\frac{1}{4}Y_1,\\
E_3&= Y_1,\\
E_4&=Y_1+Y_2,\\
E_5&=Z_1,\\
E_6&=\frac{1}{2}Z_1+Z_2,\\
E_7&=W.
\endalign
$$

Let $\Phi$ be the automorphism of $\gp$ defined on the Lie algebra level by 
$$
\align
X_1& \rightarrow -X_1+X_2+\frac{1}{4}Y_1+\frac{1}{2}Y_2,\\
X_2& \rightarrow X_2-\frac{1}{2}Y_1+\frac{1}{4}Z_1,\\
Y_1& \rightarrow -Y_1,\\
Y_2& \rightarrow 2Y_1+Y_2+Z_2,\\
Z_1& \rightarrow Z_1+\frac1{2}W,\\
Z_2& \rightarrow -Z_1-Z_2+\frac1{4}W,\\
W& \rightarrow -W.
\endalign
$$
 
A straightforward calculation shows that $\Phi_*([U,V])=[\Phi_*(U), \Phi_*(V)]$ 
for all $U,V$ in $\alg.$  Thus $\Phi$ is indeed a Lie group automorphism.

Let $\G1$ be the cocompact, discrete subgroup of $\gp$ generated by $$\{\exp(2X_1),\exp(2X_2),\exp(Y_1),\exp(Y_2),\exp(Z_1),\exp(Z_2),\exp(W)\},$$
and let $\G2=\Phi(\G1).$  Note that $\G1 \cap Z(\gp) = \G2 \cap Z(\gp) = \{\exp(jW):j\in\Z\}.$  

Let $\bar\Phi$ be the projection of $\Phi$ onto $ \qgp.$
Then $\bar\Phi$ factors as $\bar\Phi = \Psi_1 \circ \Psi_2$ where $\Psi_1$ is the 
automorphism of $\qgp$ given on the Lie algebra level by 
$$\align \BX_1 & \rightarrow -\BX_1+\BX_2+\frac{1}{4}\BY_1+\frac{1}{2}\BY_2, \\
 \BX_2 & \rightarrow \BX_2-\frac{1}{2}\BY_1,\\
 \BY_1 & \rightarrow -\BY_1,\\
 \BY_2 & \rightarrow 2\BY_1+\BY_2,\\
 \BZ_1 & \rightarrow \BZ_1,\\
 \BZ_2 & \rightarrow -\BZ_1-\BZ_2,
\endalign$$
and $\Psi_2$ is the automorphism of $\qgp$ given on the Lie algebra level by 
$$\align \BX_1 & \rightarrow \BX_1,\\
 \BX_2 & \rightarrow \BX_2+\frac{1}{4}\BZ_1,\\
 \BY_1 & \rightarrow \BY_1,\\
 \BY_2 & \rightarrow \BY_2-\BZ_1-\BZ_2,\\
 \BZ_1 & \rightarrow \BZ_1,\\
 \BZ_2 & \rightarrow \BZ_2.
\endalign$$

By rewriting $\Psi_1$ in terms of the orthonormal basis 
$\{\BE_1, \BE_2, \BE_3, \BE_4, \BE_5, \BE_6 \}$ of 
$\qmetric,$
one easily sees that 
$\Psi_1(\BE_i)=\pm\BE_i$ for 
$i=1, \dots, 6.$  
Thus the automorphism 
$\Psi_1$ is also an isometry of 
$\qGG.$ A simple calculation shows that 
$\Psi_2$ is an almost inner automorphism of 
$ \qgp.$ 
Thus by (3.2.1), $\bar\Phi$ marks the length spectrum between
$(\qnilmfld1, \qmetric)$ and $(\qnilmfld2, \qmetric).$
By (3.3.1), $\Phi$ marks the length spectrum between
$(\nilmfld1, g)$ and $(\nilmfld2, g).$

By (3.2.2), $(\nilmfld1, g)$ and $(\nilmfld2,g)$
must be isospectral on functions. 

In contrast, we
have the following.

\proclaim{4.1 Theorem \cite{Gt3, Proposition 4.11}}
The manifolds $(\nilmfld1, g)$ and 
\newline
$(\nilmfld2, g)$ are not isospectral on one-forms. 
\endproclaim
\bigskip

\subheading{\S5 The (Marked) Length Spectrum and Previous Examples}
\bigskip

We now compare the length spectra and 
marked length spectra 
of Examples I-IV described in Table I.  
The spectrum on functions, spectrum on one-forms, quasi-regular
representations and fundamental groups of these examples
were studied in \cite{Gt3}. 

We use the notation of Section 3.
\medskip

All of these examples are described by 
Theorem 3.2.3. In particular,  Examples I-IV
have the property 
$\G1 \cap Z(\gp) = \G2 \cap Z(\gp).$  

\medskip

Let $\lambda \in [L]\hy\spec (\nilmfld{i}, g ).$  
Let $m_i(\lambda)$ denote the multiplicity of $\lambda$ in 
\break
$[L]\hy\spec (\nilmfld{i},g ).$  
We decompose $m_i(\lambda)$ as 
$$m_i(\lambda) = m'_i(\lambda) + m''_i(\lambda)\tag5.1$$  
where $m''_i(\lambda)$ is the number of central 
free homotopy classes in which $\lambda$ occurs, 
and $m'_i(\lambda)$ is the number of noncentral 
free homotopy classes in which $\lambda$ occurs.

\proclaim{5.2 Proposition} For pairs of 
isospectral manifolds constructed using 
Theorem 3.2.3, the central multiplicities 
are equal; that is, 
$m''_1(\lambda) = m''_2(\lambda).$
\endproclaim

\demo{Proof of Proposition 5.2}

If $\gamma \in \G1\cap Z(\gp) = \G2\cap Z(\gp),$
then by (2.1.1),
$\lambda \in [\gamma]_\G1$ if and only if $\lambda \in [\gamma]_\G2.$
As the conjugacy classes of $\gamma$ in $\G1$ and $\G2$
respectively contain only the element $\gamma,$
we have a natural correspondence between the central
conjugacy classes in $\G1$ containing a closed geodesic of
length $\lambda$ and the central conjugacy classes in
$\G2$ containing a closed geodesic of length $\lambda.$
$\endpf$
\enddemo

Thus, for the examples below, 
we need only compare 
$m'_1(\lambda)$ and 
$m'_2(\lambda).$

\bigskip


\subheading{Example I: Remarks}
\bigskip

Let 
$$\alg =  span_{\R} \{X_1, X_2, Y_1, Y_2, Z_1, Z_2, W  \}$$
with Lie brackets
$$ [X_1,Y_1 ] = 
 [X_2, Y_2 ] = Z_1$$
$$ [X_1,Y_2 ] = Z_2$$
$$ [X_1,Z_1 ] = 
 [X_2, Z_2 ] = 
 [Y_1, Y_2 ] = W$$
and all other basis brackets zero.

Clearly $\alg$ is a strictly nonsingular, three-step
nilpotent Lie algebra.

Let
$\G1$ be the cocompact, discrete subgroup of 
$\gp$ generated by
$$ \{\exp(2X_1), \exp(2X_2), \exp(Y_1),
\exp(Y_2), \exp(Z_1), \exp(Z_2), \exp(W) \},$$
and let $\G2$ be the cocompact, discrete 
subgroup of $\gp$ generated by
$$ \{ \exp(2X_1), \exp(2X_2), \exp(Y_1),
\exp(Y_2+\frac{1}{2}Z_2), \exp(Z_1), \exp(Z_2), \exp(W)  \}.$$

The fundamental groups and the quasi-regular representations of Example I
are studied extensively in \cite{Gt1}.  There we showed that $\G1$
and $\G2$ are not abstractly isomorphic, hence $(\nilmfld1, g)$
and $(\nilmfld2, g)$ cannot possibly have the same marked length 
spectrum for any choice of left invariant metric.

Let $g$ be the left invariant metric on $\gp$ 
defined by letting 
$$ \{X_1,X_2,Y_1,Y_2,Z_1,Z_2,W \}$$ 
be an orthonormal basis of $\alg.$

In \cite{Gt2}, we showed that $(\nilmfld1, g)$ and $(\nilmfld2,g)$
do not even have the same length spectrum.  Although the same lengths
of closed geodesics occur, the multiplicities of certain lengths differ.

Example I provided the first example of a pair of representation equivalent
subgroups of a solvable Lie group producing manifolds with unequal
length spectra.  Note that nilpotent Lie groups are necessarily solvable.
\bigskip

\subheading{Example II: The (Marked) Length Spectrum}
\bigskip

Let
$$\alg = span_{\R}\{X_1,Y_1, Y_2, Z, W \}$$
with Lie brackets given by
$$[X_1,Y_1] = Z$$
$$[X_1,Z] = [Y_1, Y_2] = W$$
and all other basis brackets zero.

Clearly $\alg$ is a strictly nonsingular, three-step nilpotent Lie algebra.

Let $\G1$ be the cocompact, discrete subgroup of 
$\gp$ generated by 
$$ \{\exp(2X_1), \exp(Y_1), 
\exp(Y_2), \exp(Z), \exp(W) \}$$
and let $\G2$ be the cocompact, discrete subgroup of $\gp$ generated by
$$ \{\exp(2X_1), \exp(Y_1+\frac1{2}Z), 
\exp(Y_2), \exp(Z), \exp(W) \}.$$
Note that these generating sets are canonical in the sense that
every element of $\G1$ can be written in the form
$\exp(2n_1X_1)\exp(m_1Y_1)\exp(m_2Y_2)\exp(kZ)\exp(jW)$
for some integers $n_1, m_1, m_2, k,j.$ Likewise for $\G2.$

\proclaim{5.3 Proposition}
The above nilmanifolds have the same length spectrum, 
that is 
$$[L]\hy\spec(\nilmfld1, g) =[L]\hy\spec(\nilmfld2,g)$$ 
for any choice of left invariant metric 
$g$ of $\gp.$
\endproclaim

We showed in \cite{Gt3} that $\G1$ and $\G2$ are isomorphic as groups. 
Thus a natural question to ask is, if a pair of nilmanifolds
have the same length spectrum and have isomorphic fundamental
groups, must they necessarily have the same marked length
spectrum?  We know already from \cite{Gt2} that this need not
be true even in the two-step case.  This example is a higher-step
example with the same property.

\proclaim{5.4 Proposition}
The manifolds $(\nilmfld1,g)$ and $(\nilmfld2,g)$ 
do not have the same marked length spectrum
for any choice of left invariant metric $g$ on $\gp.$
\endproclaim

\demo{Proof of Proposition 5.4}

Let $g$ be any left invariant metric on $\gp,$ and
assume $\Psi:\G1\rightarrow\G2$ marks the length spectrum between 
$(\nilmfld1, g)$ and $(\nilmfld2,g).$
Extend $\Psi$ to the Lie group 
isomorphism 
$\Psi:\gp\rightarrow\gp$ such that 
$\Psi(\G1)=\G2.$

We showed in \cite{Gt3} Proposition 4.6
that any isomorphism 
$\Psi:\G1\rightarrow\G2$ must be given at the Lie algebra level by:

$\Psi_*(W) = \pm W,$

$\Psi_*(Z) = \pm Z + h_0W $

$\Psi_*(Y_2) = \pm Y_2$ mod $\dalg1$ 

$\Psi_*(Y_1) = \pm (Y_1 + \frac1{2}Z) + h_1 Y_2 + h_2Z$ mod 
$\dalg2$

$\Psi_*(X_1) = \pm X_1 +  \frac1{2}h_3Y_1 + \frac1{2}h_4Y_2$ mod 
$\dalg1$

\flushpar where $h_0, \ h_1,\ h_2, \ h_3$ and $h_4$ are integers
and $h_3^2+ h_4^2 \neq 0.$

By Corollary 3.2.5 and Theorem 3.2.1,
$\bar\Psi = \Phi_1 \circ \Phi_2,$ 
where $\Phi_1:\qgp\rightarrow\qgp$
is an isomorphism that is also an isometry
of $(\qgp, \qmetric),$
and $\Phi_2 \in \qG1\text{-}\AIA(\qgp).$
As $\BY_1$ and $\BY_2$ are not in 
$[\BX_1,\qalg],$
we must have

${\Phi_1}_*(\BX_1) = \pm \BX_1 + 
\frac1{2}h_3\BY_1 + \frac1{2}h_4\BY_2 + 
z_1\BZ,$ 

${\Phi_1}_*(\BY_1) = \pm \BY_1 + h_1 \BY_2 + z_2\BZ,$ 

${\Phi_1}_*(\BY_2) = \pm \BY_2 +z_3\BZ,$

${\Phi_1}_*(\BZ) = \pm \BZ,$

\flushpar for some $z_1,z_2,z_3 \in \R.$

Now $\Phi_1$ an isometry implies that for all $\BU, \BV$ in 
$\qalg,$ 
$$\left<\BU, \BV\right> = 
\left<{\Phi_1}_*(\BU),{\Phi_1}_*(\BV)\right>.\tag{$*$}$$

Letting $\BU=\BZ$ and $\BV=\BY_2$ in $(*)$, we see that
$z_3=0.$  Letting $\BU=\BY_2$ and $\BV=h_1\BY_2+z_2\BZ$
in $(*)$, 
we obtain $h_1=z_2=0.$  Finally by letting 
$\BU=\frac1{2}h_3\BY_1 + \frac1{2}h_4\BY_2 + z_1\BZ$ and 
$\BV=X_1$ in $(*)$ we see that $z_1=h_2=h_4=0,$ which contradicts
$h_3^2+ h_4^2 \neq 0.$ $\endpf$
\enddemo

Before proving Proposition 5.3, we need the following.

\proclaim{5.5 Proposition (see \cite{Gt3,Proposition 2.1}) }
Let $\G1$ and $\G2$ be cocompact, discrete subgroups of the 
Lie group $\gp$ with left invariant metric $g.$
If for each $x$ in $\gp$ we have 
$$\#\left\{\ [\gamma]_{\G1} \subset [x]_\gp\right\} = 
\#\left\{\ [\gamma]_{\G2} \subset [x]_\gp\right\},$$
then  $$[L]\hy\spec(\nilmfld1 ,g) = [L]\hy\spec(\nilmfld2 ,g).$$
Here $\#\left\{\ [\gamma]_{\G{i}} \subset [x]_G\right\}$ 
denotes the number of distinct conjugacy classes in $\G{i}$ contained in the conjugacy class of $x$ in $\gp.$
\endproclaim

\demo{Proof of Proposition 5.3}

Let $x\in\gp.$  We count the number of distinct conjugacy classes
in $\G1$ and $\G2$ contained in $[x]_\gp.$

Let 
$\gamma_1= \exp(2n_1X_1) \exp(m_1Y_1) \exp(m_2Y_2)
\exp(kZ) \exp(jW) \in \G1$ for 
\break
$ n_1, m_1, m_2, k \in \Z.$
Define the mapping $F:\G1\rightarrow\G2$ by 
$$F(\gamma_1)=\exp(2n_1X_1) \exp(m_1(Y_1+\frac1{2}Z)) 
\exp(m_2Y_2) \exp(kZ) \exp(jW).$$
The mapping $F$ gives us a correspondence between
the elements of $\G1$ and the elements of $\G2.$
Note that $F$ is {\it not} a Lie group
isomorphism.  

Now $\gamma_1$ and $F(\gamma_1)=\gamma_2$ are conjugate
in $\gp.$
In particular,
$F(\gamma_1)=a\gamma_1 a^{-1}$ where $a=e$ if $m_1=0$, and
$a=\exp(\frac{1}{2}X_1)\exp((\frac{1}{8}+\frac{k}{2m_1})Y_2)$
if $m_1\neq0.$  
Thus $[\gamma_1]_\G1 \subset [x]_\gp$
if and only if $[F(\gamma_1)]_\G2 \subset [x]_\gp.$

To use Proposition 5.5, we must now compare the number
of distinct conjugacy classes in $\G1$ and $\G2$ respectively
that are contained in a fixed $[x]_\gp.$

Using the Camp\-bell-Baker-Haus\-dorff
formula, two elements 
$$\gamma_1 = 
\exp(2n_1X_1) \exp(m_1Y_1) \exp(m_2Y_2)
\exp(kZ) \exp(jW)$$
and $$\gamma'_1=\exp(2n'_1X_1) \exp(m'_1Y_1) \exp(m'_2Y_2)
\exp(k'Z) \exp(j'W),$$ of $\G1$
are conjugate in $\G1$ if and only if there exist
integers $\b{n}{1}, \b{m}{1}, \b{m}{2}, \b{k}{ }$ such that
$$\split 
n'_1 &= n_1, \quad m'_1 = m_1, \quad m'_2 = m_2,\\
k' &= k + 2 m_1 \bar{n}_1 - 2 n_1 \bar{m}_1,\\
j' &= j + m_2\bar{m}_1 - m_1\bar{m}_2 + 
2k\bar{n}_1 - 2n_1\bar{k} \\
&\qquad + 2m_1{\bar{n}_1}^2 - 4n_1\bar{n}_1\bar{m}_1 + 
2{n_1}^2\bar{m}_1.\endsplit$$

Let $K = gcd(2n_1,2m_1).$ 
From the above, we see that every conjugacy class in $\G1$  
contains at least one representative such that 
$k \in \{ 1, 2, \cdots, K\}.$
We call such a representative {\it nice.}
Two nice representatives
are in the same conjugacy class in 
$\G1$ if and only if $k = k'$ and
there exist integers $\b{n}{1}, \b{m}{1}, \b{m}{2}, \bar{k}$
such that  $m_1 \bar{n}_1 - n_1\bar{m}_1 = 0$
and
$$\split j' &= j + m_2\bar{m}_1 - m_1\bar{m}_2 + 
2k\bar{n}_1 - 2n_1\bar{k} \\
&\qquad + 2m_1{\bar{n}_1}^2 - 4n_1\bar{n}_1\bar{m}_1 + 
2{n_1}^2\bar{m}_1 \endsplit$$

Similarly, two elements of $\G2$
$$\gamma_2 = \exp(2n_1X_1) \exp(m_1(Y_1+\frac1{2}Z)) 
\exp(m_2Y_2) \exp(kZ) \exp(jW),$$ and 
$$\gamma'_2 =  \exp(2n'_1X_1) \exp(m'_1(Y_1+\frac1{2}Z)) 
\exp(m'_2Y_2) \exp(k'Z) \exp(j'W),$$
are conjugate in $\G2$ if and only if there exist integers
$\b{n}{1}, \b{m}{1}, \b{m}{2}, \b{k}{ }$ so that 
$$\split 
n'_1 &= n_1, \quad m'_1 = m_1, \quad m'_2 = m_2,\\
k' &= k + 2 m_1 \bar{n}_1 - 2n_1\bar{m}_1,\\
j' &= j + (m_1\bar{n}_1 - n_1\bar{m}_1)+ 
m_2\bar{m}_1 - m_1\bar{m}_2 + 2k\bar{n}_1 - 2n_1\bar{k} \\
&\qquad + 2m_1{\bar{n}_1}^2 - 4n_1\bar{n}_1\bar{m}_1 + 
2{n_1}^2\bar{m}_1.\endsplit$$

Again we see that  every conjugacy class in $\G2$ contains  
at least one nice representative, that is,
a representative such that 
$k \in \{ 1, 2, \cdots, K\},$ where 
$K=gcd(2n_1,2m_1)$ as above.  
Again, two nice representatives
are in the same conjugacy class in 
$\G2$ if and only if $k = k'$ 
and 
there exist integers $\b{n}{1}, \b{m}{1}. \b{m}{2}, \bar{k}$
such that $m_1 \bar{n}_1 - n_1\bar{m}_1 = 0$ and  
$$\split j' &= j + m_2\bar{m}_1 - m_1\bar{m}_2 +
 2k\bar{n}_1 - 2n_1\bar{k} \\
&\qquad + 2m_1{\bar{n}_1}^2 - 4n_1\bar{n}_1\bar{m}_1 + 
2{n_1}^2\bar{m}_1 \endsplit$$

Note that the correspondence $ F:\G1\rightarrow\G2$
sends nice representatives to nice representatives. 
Thus if we restrict ourselves to nice representatives,
the conjugacy conditions are equivalent.
That is, two nice representatives $\gamma_1$ and $\gamma'_1$
are in the same conjugacy class in $\G1$ if and only if
the corresponding elements $ F(\gamma_1)$
and $ F(\gamma_2)$ are in the same
conjugacy class in $\G2.$  

Let $\gamma_1, \gamma_2, \dots, \gamma_L$ be nice
representatives of the $L$ distinct conjugacy classes
in $\G1$ contained in $[x]_\gp.$
Then $ F(\gamma_1),  F(\gamma_2), \dots,  F(\gamma_L)$
are nice representatives of $L$ distinct conjugacy classes in $\G2.$
The same applies to $ F^{-1}:\G2\rightarrow\G1.$

Thus 
$$\#\left\{\ [\gamma]_{\G1} \subset [x]_\gp\right\} = 
\#\left\{\ [\gamma]_{\G2} \subset [x]_\gp\right\},$$
as desired.\endpf
\enddemo
\bigskip


\subheading{Example III:  The Length Spectrum}
\bigskip

Let 
$$\alg =  span_{\R} \{X_1, X_2, Y_1, Y_2, Z_1, Z_2, W  \}$$
with Lie brackets
$$ [X_1,Y_1 ] = 
 [X_2, Y_2 ] = Z_1$$
$$ [X_1,Y_2 ] = Z_2$$
$$ [X_1,Z_1 ] = 
 [X_2, Z_2 ] = 
 [Y_1, Y_2 ] = W$$
and all other basis brackets zero.

Clearly $\alg$ is a strictly nonsingular, three-step
nilpotent Lie algebra.

Let
$\G1$ be the cocompact, discrete subgroup of 
$\gp$ generated canonically by
$$ \{\exp(2X_1), \exp(2X_2), \exp(Y_1),
\exp(Y_2), \exp(Z_1), \exp(Z_2), \exp(W) \},$$
and let $\G2$ be the cocompact, discrete 
subgroup of $\gp$ generated canonically by
$$ \{ \exp(X_1), \exp(X_2), \exp(2Y_1),
\exp(2Y_2), \exp(Z_1), \exp(Z_2), \exp(W)  \}.$$

Let $g$ be the left invariant metric on $\gp$ 
defined by letting 
$$ \{X_1,X_2,Y_1,Y_2,Z_1,Z_2,W \}$$ 
be an orthonormal basis of $\alg.$

\proclaim{5.6 Proposition }
The nilmanifolds $(\nilmfld1, g)$ and 
$(\nilmfld2,g)$ do not have the same length spectrum.
In particular,  the multiplicity of the length 1 in 
\break
$[L]\hy\spec(\nilmfld1,g)$ 
is greater than its multiplicity in
$[L]\hy\spec(\nilmfld2,g).$
\endproclaim

\demo{Proof of Proposition 5.6}

By Proposition 5.2 we need only consider
the non\-central free ho\-mo\-to\-py 
\break
classes. That is,
we need only show $m'_1(1) >  m'_2(1).$  

Let 
$$\eightpoint{\gamma = \exp(A_1n_1X_1) \exp(A_2n_2X_2) 
\exp(B_1m_1Y_1) \exp(B_2m_2Y_2)
\exp(k_1Z_1) \exp(k_2Z_2) \exp(jW)}$$
for integers 
$n_1,n_2,m_1,m_2,k_1,k_2,j$ and $A_1, A_2, B_1, B_2 \in \{1,2\}.$ 
Note that $\gamma \in \G1$  if and only if
$$A_1=A_2=2,\quad  B_1=B_2=1\tag{$*$}$$ 
and if $\gamma \in \G2$  if and only if
$$A_1=A_2=1, \quad B_1=B_2=2.\tag{$**$}$$  

By Theorem 3.2.4, to determine if 
$1\in [\gamma]_\G{i}$ for noncentral $\gamma\in\G{i},$
we need only determine if 
$1\in [\bar\gamma]_\qG{i}.$  
That is, rather than  looking at the lengths of 
closed geodesics on the three-step nilmanifolds 
$(\nilmfld{i},g),$  we instead look at 
the lengths of closed geodesics on the 
quotient two-step nilmanifolds 
$( \qG{i}\backslash \qgp,\qmetric)$ for $i=1,2.$  

The Lie algebra of $\qgp$ is 
$\qalg = \alg/\dalg2=
 span_\R\{\BX_1, \BX_2, \BY_1, \BY_2, \BZ_1, \BZ_2\}$
with Lie brackets
$$ [\BX_1,\BY_1 ] = 
 [\BX_2, \BY_2 ] = \BZ_1$$
$$ [\BX_1,\BY_2 ] = \BZ_2,$$
and all other basis brackets zero.

We may now  use the following result due to Eberlein.  

\proclaim{5.7 Theorem \cite{E1}}
Let $N$ be a simply connected, two-step 
nilpotent Lie group with Lie algebra $\nn$
and left invariant metric $g.$
Let $\GG$ be a cocompact, discrete subgroup of $N.$
Let $\center$ be the center of $\nn$ and $\vv$ 
the orthogonal complement of $\center$ in 
$\nn.$  Any element $\gamma \in \GG$ may be expressed 
uniquely as $\exp(V^*+Z^*)$ where $V^*\in\vv$ 
and $Z^*\in\center.$  Let $Z^{**}$ be the 
component of $Z^*$ orthogonal to $[V^*,\nn].$  
Let $\lambda>0.$

(1) If $\lambda\in[\gamma]_\GG$, then $|V^*| \leq \lambda 
\leq \sqrt{|V^*|^2 + |Z^{**}|^2}.$ 

(2) The period $\lambda=|V^*| \in [\gamma]_\GG$ 
if and only if $|Z^{**}|=0.$

(3) The period $\lambda=\sqrt{|V^*|^2+|Z^{**}|^2} \in [\gamma]_\GG.$
\endproclaim 

Here $\nn = \qalg$ and the metric 
$\qmetric$ is determined by the orthonormal basis  
of $\qalg$ 
$$ \{\BX_1,\BX_2,\BY_1,\BY_2,\BZ_1,\BZ_2 \}.$$

By Theorem 5.7, to find $\bar\gamma$ such that $1\in [\bar\gamma]_\qG{i},$ we need 
$\bar\gamma = \exp( \BV^* + \BZ^*)$ 
such that 
$| \BV^*|^2 \leq 1 \leq  | \BV^*|^2 + |\BZ^{**}|^2,$ 
where $ \BV^* \in  span_\R\{\BX_1,\BX_2,\BY_1,\BY_2\}$ and 
$\BZ^* \in  span_\R\{\BZ_1,\BZ_2\}.$

For both $\G1$ and $\G2,$  
$ \BV^* = A_1n_1\BX_1 + A_2n_2\BX_2 + B_1m_1\BY_1 + B_2m_2\BY_2, $ 
where 
$n_1,n_2,m_1,m_2 \in \Z$ 
Note that if $|V^*|\neq0,$
$| \BV^*|^2 = A_1^2n_1^2 + A_2^2n_2^2 + B_1^2m_1^2 + B_2^2m_2^2 \geq1.$  
So $| \BV^*|^2 \leq 1$ if and only if $| \BV^*|^2 =1.$
By Theorem 5.4,  $\lambda = 1= |\BV^*| \in [\bar\gamma]_\qG{i}$ 
if and only if $|\BZ^{**}|=0.$
\medskip

So if $\bar\gamma = \exp( \BV^* + \BZ^*)$ with 
$|\BV^*| \not=0,$ 
then $1\in  [\bar\gamma]_\qG{i}$ if and only if 
where $|\BV^*|=1$ and $|\BZ^{**}|=0.$

We consider two cases.
\bigskip

Case 1:  $(n_1)^2+(m_2)^2 \neq 0.$
\medskip

In this case, $\bar\center =  [\log\bar\gamma, \qalg ],$ 
so $\BZ^{**}$ is automatically zero. 
Applying the condition $|\BV^*|=1$ and 
lifting to the three-step level, we have 
$1\in[\gamma_1]_\G1$ 
if and only if (see ($*$))
$$\gamma_1 = \exp(\pm Y_2) \exp(k_1Z_1) \exp(k_2Z_2) \exp(jW),$$
and   $1\in[\gamma_2]_\G2$ if and only if (see ($**$))
$$\gamma_2= \exp(\pm X_1) \exp(k_1Z_1) \exp(k_2Z_2) \exp(jW).$$

We must now compare the number of distinct free homotopy 
classes of $\G1$ and $\G2$ that take on one of these forms.

Another element
$\gamma'_1=\exp(\pm Y_2)\exp(k'_1Z_1)\exp(k'_2Z_2)exp(j'W)$
of $\G1$ 
is conjugate to $\gamma_1$ in $\G1$ 
if and only if there exist integers 
$\b{n}{1}, \b{n}{2}, \b{m}{1}$ and $\b{k}{1}$ such that 
$$k'_1 = k_1 \pm 2\b{n}{2}; \quad k'_2 = k_2 \pm 2\b{n}{1}; \quad
j' = j \pm \b{m}{1} + 2k_1\b{n}{1} + 
2k_2\b{n}{2} \pm 4 \b{n}{1}\b{n}{2}.$$

Another element
$\gamma'_2 = \exp(\pm X_1) \exp(k'_1Z_1) \exp(k'_2Z_2) \exp(j'W)$ 
 of $\G2$ is conjugate to $\gamma_2$ in $\G2$ 
if and only if there exist integers 
$\b{n}{1}, \b{n}{2}, \b{m}{1}$ and $\b{m}{2}$ 
such that 
$$k'_1 = k_1 \mp 2\b{m}{1}; \quad k'_2 = k_2 \mp 2\b{m}{2}; \quad 
j' = j \mp \b{k}{1} + \b{m}{1}+ k_1\b{n}{1} + 
k_2\b{n}{2} \mp 2 \b{m}{1}\b{n}{1} \mp 2 \b{m}{2}\b{n}{2}.$$

For $\G1$ we have two choices 
$\{-1, +1\}$ for the coefficient of $Y_2,$ 
two choices for $k_1,$  
two choices for $k_2$ and 
one choice for $j$ for a total of 
$8$ distinct free homotopy classes.
For $\G2$ we have two choices 
$\{-1, +1\}$ for the coefficient of $X_1,$ 
two choices for $k_1,$  
two choices for $k_2$ 
and one choice for $j$ 
for a total of $8$ distinct free homotopy classes.
Thus, the multiplicities of 1 coming from this case are equal.
\bigskip

Case 2: ${n_1}^2 + {m_2}^2 = 0$ but ${n_2}^2 + {m_1}^2 \neq 0.$
\medskip

In this case, $[\log\bar\gamma, \qalg ]= span_\R\{\BZ_1\},$ 
so $\BZ^{**}=0$ if and only if $k_2=0.$
Applying the condition $|\BV^*|=1$ and
lifting to the three-step level, we have 
$1\in[\gamma_1]_\G1$ 
if and only if (see ($*$))
$$\gamma_1 = \exp(\pm Y_1) \exp(k_1Z_1) \exp(jW),$$
and  $1\in[\gamma_2]_\G2$ if and only if (see ($**$))
$$\gamma_2 = \exp(\pm X_2) \exp(k_1Z_1) \exp(jW).$$

We must now count the number of distinct free homotopy classes of 
$\G1$ and $\G2$ that take on one of these forms.

Another element 
$\gamma'_1= \exp(\pm Y_1) \exp(k'_1Z_1) \exp(j'W)$ 
of $\G1$ is conjugate to 
$\gamma_1$ in $\G1$ 
if and only if there exist integers 
$\b{n}{1}, \b{m}{2}$ such that 
$$k'_1 = k_1 \pm 2\b{n}{1}; \qquad 
j' = j \mp \b{m}{2} + 2k_1\b{n}{1} \pm 2{\b{n}{1}}^2.$$

Another element $\gamma'_2  = \exp(\pm X_2) \exp(k'_1Z_1) \exp(j'W)$ 
in $\G2$ is conjugate to $\gamma_2$ in $\G2$ 
if and only if there exist integers 
$\b{n}{1}, \b{m}{2}$ and $\b{k}{2}$ such that 
$$k'_1 = k_1 \mp 2\b{m}{2}; \qquad 
j' = j \mp \b{k}{2} + 
k_1\b{n}{1} \mp 2{\b{n}{1}}\b{m}{2}$$

For $\G1$ we have two choices $\{-1, +1\}$ 
for the coefficient of $Y_1,$ 
two choices for $k_1,$ 
and one choice for $j$ 
for a total of $4$ distinct free homotopy classes.
For $\G2$  we have two choices $\{-1, +1\}$ 
for the coefficient of $X_2,$ 
two choices for $k_1,$ 
and one choice for $j$ 
for a total of $4$ distinct free homotopy classes.
Again, the multiplicities of 1 coming from this case are equal.
\bigskip

Case 3: $|\BV^*|=0, |\BZ^*|\not=0.$
\medskip

Let  $\gamma = \exp(k_1Z_1) \exp(k_2Z_2) \exp(jW),$
for $k_1, k_2, j\in\Z.$ 
Note that $\gamma \in \G1 \cap \G2.$ 
Thus by (2.1.1), any period occurring in $[\gamma]_\G1$ will also 
occur in $[\gamma]_\G2.$
Let $\gamma' =  \exp(k'_1Z_1) \exp(k'_2) \exp(j'W)$
be another element of $\G1\cap\G2,$ where $k'_1, k'_2, j' \in \Z.$

Now $\gamma'$ is conjugate to $\gamma$ in $\G1$ 
if and only if there exists integers 
$\b{n}{1}, \b{n}{2}$ such that 
$$k'_1 = k_1; \qquad k'_2 = k_2; \qquad 
j' = j + 2(k_1\b{n}{1} + k_2\b{n}{2}).$$

However
$\gamma'$ is conjugate to $\gamma$ in $\G2$ 
if and only if there exists integers 
$\b{n}{1}, \b{n}{2}$ such that 
$$k'_1 = k_1; \qquad k'_2 = k_2; \qquad 
j' = j + (k_1\b{n}{1} + k_2\b{n}{2}).$$

Note that there are twice as many distinct conjugacy
classes represented by elements of the form 
$\gamma=\exp(k_1Z_1) \exp(k_2Z_2) \exp(jW)$ for $\G1$ as for $\G2.$
Thus to show the multiplicities are not equal here, 
we need to exhibit a closed geodesic of length 1 
in just one free homotopy class of this form.

Note that $|\BZ^{**}|^2=|\BZ^*|^2=k_1^2+k_2^2.$
By Theorem 5.7(iii) and lifting to the three-step
level, we see $\sqrt{k_1^2+k_2^2} \in [\gamma]_\G1$
and $\sqrt{k_1^2+k_2^2} \in [\gamma]_\G2.$
Thus
$1 \in [\exp(\pm Z_i)]_\G{j},$ $i,j=1,2.$
 
Therefore, for Case 3,   
$1$ occurs with twice the multiplicity in 
$[L]\hy\spec(\nilmfld1,g )$ as it does in 
$[L]\hy\spec(\nilmfld2,g ).$

As the multiplicities of 1 are equal in all of the other cases,
the multiplicities of 1 is {\it not} equal, as claimed. $\endpf$
\enddemo

\bigskip


\subheading{Example IV:  The Length Spectrum}
\bigskip

Here the Lie algebra is the same Lie algebra as
Example II, that is
$$\alg = span_{\R} \{X_1, Y_1, Y_2, Z, W  \}$$
with Lie brackets
$$ [X_1,Y_1 ] = Z$$
$$ [X_1,Z ] =  [Y_1, Y_2 ] = W$$
and all other basis brackets zero.

Let
$\G1$ be the cocompact, discrete subgroup of 
$\gp$ generated canonically by
$$ \{ \exp(2X_1), \exp(Y_1), \exp(Y_2),
\exp(Z), \exp(W) \},$$
and let $\G2$ be the cocompact, 
discrete subgroup of $\gp$ 
generated canonically by
$$ \{\exp(X_1), \exp(2Y_1), \exp(2Y_2), 
\exp(Z), \exp(W)  \}.$$

Let $g$ be the left invariant metric on 
$\gp$ defined by letting 
$$ \{X_1,Y_1,Y_2,Z,W \}$$ 
be an orthonormal basis of $\alg.$

\proclaim{5.8 Proposition}
The nilmanifolds $ (\nilmfld1, g )$ and 
$ (\nilmfld2,g )$ do not have 
the same length spectrum.
In particular, the multiplicity of the length 
\break
$\lambda = \sqrt{4\pi(7-\pi)}$
in 
$[L]\hy\spec(\nilmfld1,g)$ 
is greater than its multiplicity in 
\break
$[L]\hy\spec(\nilmfld2,g).$
\endproclaim

\demo{Proof of Proposition 5.8}

By Proposition 5.2, we only consider 
the noncentral 
free homotopy classes. In particular,  we show 
$m'_1(\lambda) >  m'_2(\lambda)$
where $\lambda = \sqrt{4\pi(7-\pi)}.$  

By Theorem 3.2.4 if we wish to determine if 
$\lambda \in [\gamma]_\G{i}$ for noncentral $\gamma\in\G{i},$
we need only determine if 
$\lambda \in [\bar\gamma]_\qG{i}.$  
That is, rather than  looking at the lengths of 
closed geodesics on the three-step nilmanifolds 
$(\nilmfld{i},g),$  we instead look at 
the lengths of closed geodesics on the 
quotient two-step nilmanifolds 
$( \qG{i}\backslash \qgp,\qmetric)$ for $i=1,2.$ 

However, for this example, 
$\qalg \cong \hh_1 \oplus \R$ 
where $\hh_1$ 
denotes the three-dimensional Heisenberg algebra.  
To see this, note that 
$$\hh_1 \cong  \{\BX_1,\BY_1,\BZ \}, \qquad \text{ and } \qquad [\BX_1,\BY_1]=\BZ.$$ 
This is an ideal in $\qalg.$ And 
$$\R \cong \{\BY_2\}$$ which is also an ideal in $\qalg.$ 
Let $H_1$ be the three-dimensional Heisenberg group.  
Note that 
$$H_1 \cong  \{ \exp(x_1\BX_1) \exp(y_1\BY_1) \exp(z\BZ):
x_1,y_1,z \in \R \}.$$

This direct sum is actually a {\it Riemannian} 
direct sum, as the metric may also be written as 
$$\qmetric = \qmetric_1 \oplus \qmetric_2$$ 
where $\qmetric_1$ is the left invariant metric on 
$\hh_1$ given by the orthonormal basis 
\break 
$ \{\BX_1, \BY_1,\BZ  \}$ and $\qmetric_2$ 
is the left invariant metric on $\R$ 
given by the unit vector $\{\BY_2\}.$

Furthermore, as 
$ \qG{i} =  ( \qG{i}\cap H_1 ) 
\oplus  ( \qG{i}\cap\R ),$ 
we also have the Riemannian direct sum 
$$ ( \qnilmfld{i},\qmetric ) \cong 
 ( ( \qG{i}\cap H_1 )\backslash H_1, 
\qmetric_1 ) 
\oplus 
 ( ( \qG{i}\cap\R )\backslash\R, 
\qmetric_2 ).$$

Using rescaling of geodesics, it is not difficult to show that 
$\lambda \in [\bar\gamma]_ \qG{i}$ 
if and only if 
$$\lambda^2 = \lambda_1^2 + \lambda_2^2$$ 
where $\lambda_1 \in [\bar\gamma_1]_{ \qG{i}\cap H_1}$ 
and $\lambda_2 \in [\bar\gamma_2]_{ \qG{i}\cap\R}.$
Here 
$\bar\gamma = (\bar\gamma_1, \bar\gamma_2)$ 
with respect to the direct product 
$ \qG{i} =  ( \qG{i}\cap H_1 ) 
\oplus  ( \qG{i}\cap\R ).$

Now, the length spectrum of 
$ ( ( \qG{i}\cap\R )\backslash\R, 
\qmetric_2 )$ 
is easily seen to be 
$ |\log(\bar\gamma_2) |$ 
for all 
$\bar\gamma_2 \in  \qG{i}\cap\R.$  
Thus the length spectrum here 
(not counting multiplicities) is precisely the positive integers.

The length spectrum of 
$ ( ( \qG{i}\cap H_1 )\backslash H_1, 
\qmetric_1 )$ 
has been calculated by both Gordon and Eberlein 
(see \cite {E}, \cite {G1})
and is known to be
 
\medskip

\flushpar
(i) $ |\log(\bar\gamma_1) |$ 
if 
$\bar\gamma_1 \in  \qG{i}\cap H_1,$ 
for $\bar\gamma_1 \not \in Z(H_1).$
\smallskip

\flushpar
(ii) 
$ \{  |log(\bar\gamma_1) |, 
\sqrt{ (4\pi k )
 ( |\log (\bar\gamma_1 ) |-\pi k )} 
: 1 \leq k < 
 (\frac1{2\pi}
 |\log (\bar\gamma_1 ) | ),
k\in\Z \},$ 
for  \newline $\bar\gamma_2\in\qG{i}\cap Z(H_1).$

Nonintegral lengths occur in 
$ ( \qG{i}\cap H_1 )\backslash H_1$  
only when $ |\log(\bar\gamma_1) | \geq 2\pi>6.$ 

Also note that 
$\sqrt{4\pi(7-\pi)} \in [\bar\gamma_1]_{ \qG{i}\cap H_1}$ 
if and only if  
$\bar\gamma_1 = exp(\pm 7\BZ) \in  \qG{i} \cap H_1.$ 
This is the smallest possible nonintegral length.  

Thus $$4\pi(7-\pi) = \lambda^2 = \lambda_1^2 + \lambda_2^2$$ 
if and only if 
$\lambda_2^2 = 0$ and 
$\lambda_1^2 = 4\pi(7-\pi)$
if and only if   
$\bar\gamma = exp(\pm 7\BZ) \in  \qG{i}.$

By lifting to $(\nilmfld{i},g),$  we see
$\sqrt{4\pi(7-\pi)} \in [\gamma]_\G{i}$ 
if and only if 
$$\gamma = exp(\pm 7Z)exp(jW) \in \G{i}.$$

We now count the number of 
distinct free homotopy classes represented by a 
$\gamma$ of this form.

Let $\gamma' = exp(\pm 7 Z)exp(j'W).$  

Now $\gamma'$ is conjugate to 
$\gamma$ in $\G1$ 
if and only if there exists integer 
$\b{n}{1}$ such that 
$$j' = j \pm 14\b{n}{1}.$$

However, 
$\gamma'$ is conjugate to 
$\gamma$ in $\G2$ 
if and only if there exists integer 
$\b{n}{1}$ such that 
$$j'= j \pm 7\b{n}{1}.$$

Thus there are 14 choices for 
$j$ in $\G1$ and there are 7 choices for 
$j$ in $\G2.$
So the multiplicity of the length 
$\sqrt{4\pi(7-\pi)}$ in 
$ (\nilmfld1,g )$ is 28, 
(14 for each of 
$exp(+7Z)exp(jW)$ and 
$exp(-7Z)exp(jW)$), 
and likewise the multiplicity in 
$ (\nilmfld2,g )$ is 14.

Thus the multiplicities of 
$\sqrt{4\pi(7-\pi)}$  
are {\it not} equal here, as claimed. 
$\endpf$
\enddemo

\enddocument